\documentclass[twoside]{article} 
\usepackage{graphicx}
\title{\bf On some F\'ejer-type trigonometric sums}
\author{\sc R. B. Paris\footnote{E-mail address:\ \ {\tt r.paris@abertay.ac.uk}}\\
\\
{\em Division of Computing and Mathematics,}\\
{\em Abertay University, Dundee DD1 1HG, UK}\\
}

\begin{document}
\newcommand{\bee}{\begin{equation}}
\newcommand{\ee}{\end{equation}}
\def\f#1#2{\mbox{${\textstyle \frac{#1}{#2}}$}}
\def\dfrac#1#2{\displaystyle{\frac{#1}{#2}}}
\newcommand{\fr}{\frac{1}{2}}
\newcommand{\fs}{\f{1}{2}}
\newcommand{\g}{\Gamma}
\newcommand{\al}{\alpha}
\newcommand{\bb}{\beta}
\newcommand{\om}{\omega}
\newcommand{\br}{\biggr}
\newcommand{\bl}{\biggl}
\newcommand{\ra}{\rightarrow}
\renewcommand{\topfraction}{0.9}
\renewcommand{\bottomfraction}{0.9}
\renewcommand{\textfraction}{0.05}
\newcommand{\mcol}{\multicolumn}
\newcommand{\gtwid}{\raisebox{-.8ex}{\mbox{$\stackrel{\textstyle >}{\sim}$}}}
\newcommand{\ltwid}{\raisebox{-.8ex}{\mbox{$\stackrel{\textstyle <}{\sim}$}}}
\date{}
\maketitle
\pagestyle{myheadings}
\markboth{\hfill {\it R.B. Paris} \hfill}
{\hfill {\it Asymptotics of a modified exponential integral} \hfill}
\begin{abstract} 
We examine the four F\'ejer-type trigonometric sums of the form
\[S_n(x)=\sum_{k=1}^n \frac{f(g(kx))}{k}\qquad (0<x<\pi)\]
where $f(x)$, $g(x)$ are chosen to be either $\sin x$ or $\cos x$. The analysis of the sums with $f(x)=g(x)=\cos x$, $f(x)=\cos x$, $g(x)=\sin x$ and $f(x)=\sin x$, $g(x)=\cos x$ is reasonably straightforward. It is shown that these sums exhibit unbounded growth as $n\to\infty$ and also present `spikes' in their graphs at certain $x$ values for which we give an explanation.

The main effort is devoted to the case $f(x)=g(x)=\sin x$, where we present arguments that strongly support the conjecture made by H. Alzer that $S_n(x)>0$ in $0<x<\pi$. The graph of the sum in this case presents a jump in the neighbourhood of $x=\f{2}{3}\pi$. This jump is explained and is quantitatively estimated when $n\to\infty$.

\vspace{0.4cm}

\vspace{0.3cm}

\end{abstract}

\vspace{0.2cm}

\noindent $\,$\hrulefill $\,$

\vspace{0.2cm}

\begin{center}
{\bf 1. \  Introduction}
\end{center}
\setcounter{section}{1}
\setcounter{equation}{0}
\renewcommand{\theequation}{\arabic{section}.\arabic{equation}}
In 1910 F\'ejer conjectured that all partial sums of the series
$\sum_{k=1}^\infty \frac{\sin kx}{k}=\fs(\pi-x)$ in the interval $0<x\leq\pi$
are positive, namely
\bee\label{e11}
s_n(x)=\sum_{k=1}^n\frac{\sin kx}{k} >0 \qquad n\geq 1,\ 0<x<\pi.
\ee
Proofs of this result were given by Jackson \cite{Jack}, Gronwall \cite{Gron}, Landau \cite{Lan}, F\'ejer \cite{Fej} and Tur\'an \cite{Tur}. A similar result for the cosine series is \cite{WHY}
\bee\label{e12}
c_n(x)=\sum_{k=1}^n\frac{\cos kx}{k} >-1\qquad n\geq 1,\ 0<x<\pi.
\ee

Extensions of the above trigonometric sums considered here have the form
\[S_n^{(1)}(x)=\sum_{k=1}^n\frac{\cos (\cos kx)}{k},\qquad S_n^{(2)}(x)=\sum_{k=1}^n\frac{\cos (\sin kx)}{k},\]
\bee\label{e13}
S_n^{(3)}(x)=\sum_{k=1}^n\frac{\sin (\cos kx)}{k},\qquad S_n^{(4)}(x)=\sum_{k=1}^n\frac{\sin (\sin kx)}{k}.
\ee
Recently, Alzer \cite{Alz}  made the following conjecture concerning $S_n^{(4)}(x)$ based on numerical evidence:
\vspace{0.05cm}

\noindent{\bf Conjecture}.\ {\it For all positive integer $n$ and $0<x<\pi$ we have
\bee\label{e14}
S_n^{(4)}(x)=\sum_{k=1}^n\frac{\sin (\sin kx)}{k}>0,\ee
where the constant lower bound is sharp.}
\vspace{0.2cm}

\noindent The investigation of the sums $S_n^{(j)}(x)$, $1\leq j\leq 3$ is relatively straightforward and will be dealt with in Section 3. Our main effort is concentrated on the sum $S_n^{(4)}(x)$ and is considered in Section 4. We do not claim to establish a proof of the above conjecture but we present arguments that strongly suggest its validity. 
For this purpose, we shall make use of the identities \cite[p.~226]{DLMF}
\bee\label{e15}
\left.\begin{array}{l}
\sin(z\cos \theta)={\displaystyle 2\sum_{r=0}^\infty (-)^r J_{2r+1}(z) \cos ((2r+1)\theta)},\\
\\
\sin(z\sin \theta)={\displaystyle 2\sum_{r=0}^\infty J_{2r+1}(z) \sin ((2r+1)\theta)}
\end{array}\right\}
\ee
valid for $z$, $\theta\in{\bf C}$, where $J_\nu(z)$ denotes the Bessel function of the first kind, together with some basic properties of the sine and cosine series $s_n(x)$ and $c_n(x)$ in (\ref{e11}) and (\ref{e12}), which we present in the next section.
\vspace{0.6cm}

\begin{center}
{\bf 2. \ Properties of the basic series $s_n(x)$ and $c_n(x)$}
\end{center}
\setcounter{section}{2}
\setcounter{equation}{0}
\renewcommand{\theequation}{\arabic{section}.\arabic{equation}}
In this section we present some basic properties of the simple sine and cosine series $s_n(x)$ and $c_n(x)$ defined in (\ref{e11}) and (\ref{e12}).
Both these series are periodic functions of $x$ with period $2\pi$. Typical plots of these functions are illustrated in Fig.~1. 
\begin{figure}[th]
	\begin{center}	{\tiny($a$)}\includegraphics[width=0.375\textwidth]{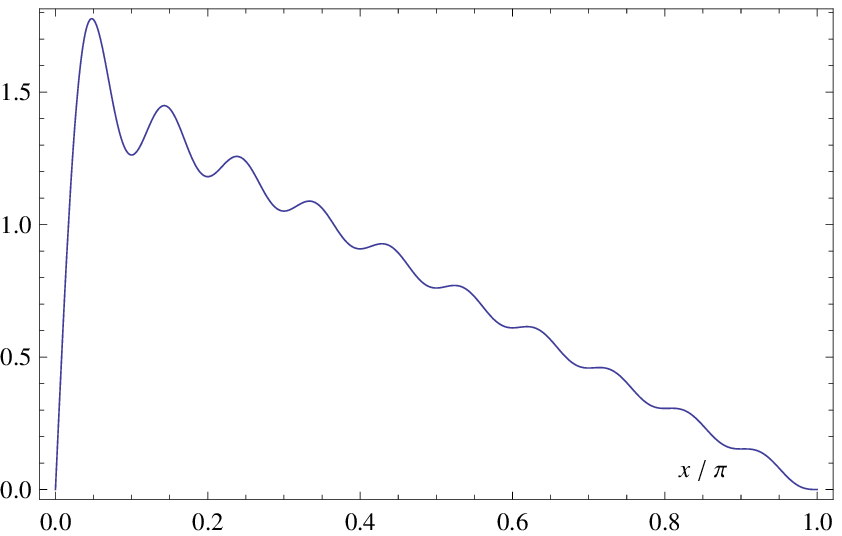}
	\qquad
	{\tiny($b$)}\includegraphics[width=0.375\textwidth]{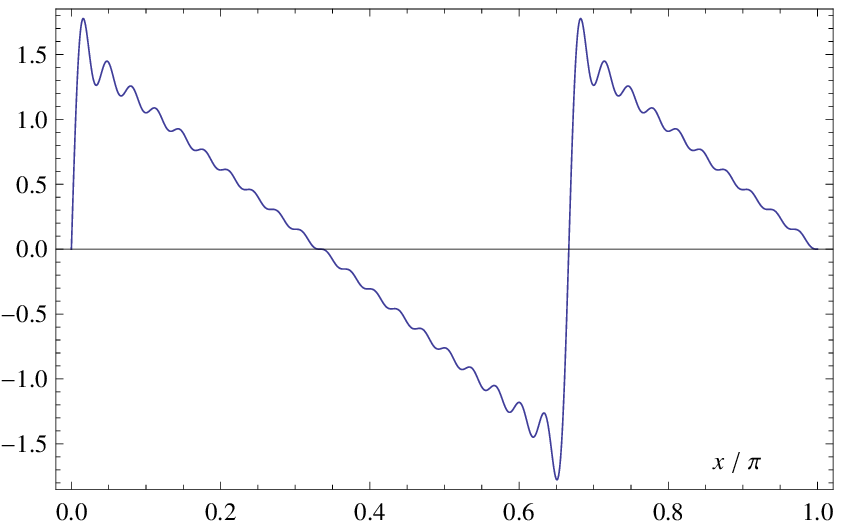} 
	\bigskip
	
		{\tiny($c$)}\includegraphics[width=0.375\textwidth]{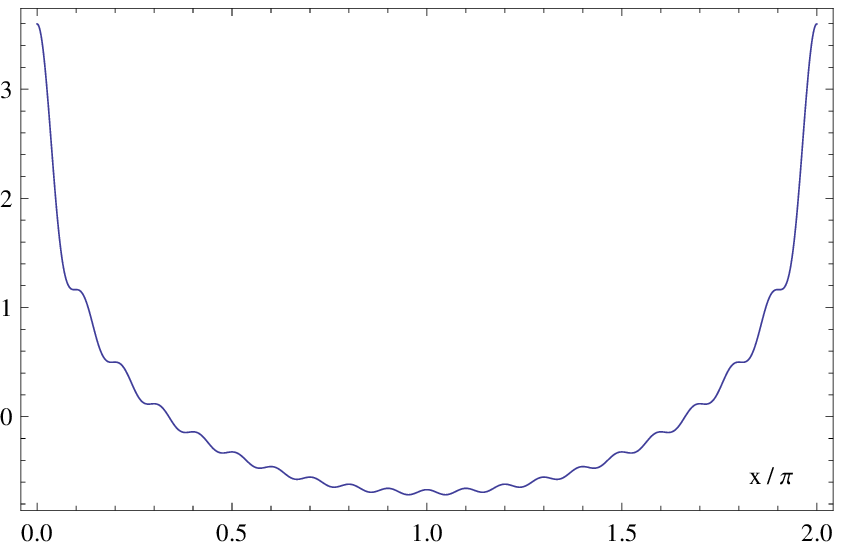}
	\qquad
	{\tiny($d$)}\includegraphics[width=0.375\textwidth]{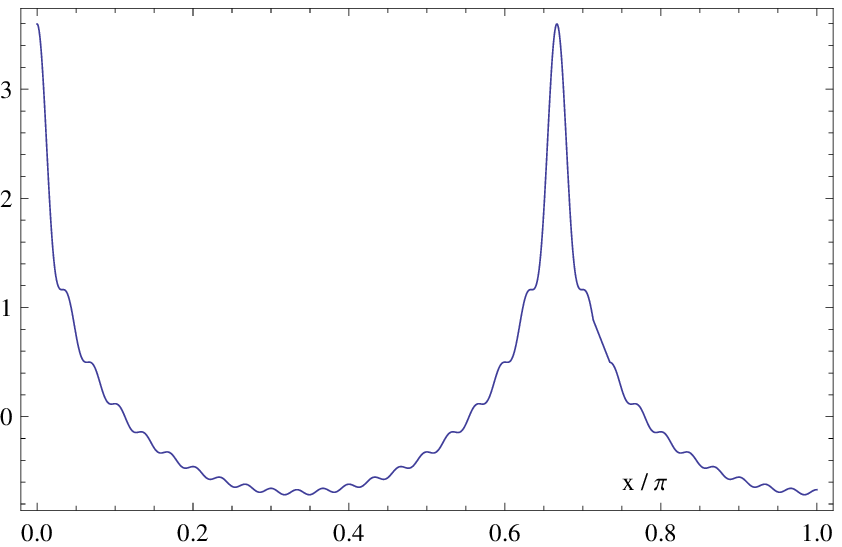} \\
	
\caption{\small{Plots of (a) $s_n(x)$, (b) $s_n(3x)$ for $0\leq x\leq\pi$, (c) $c_n(x)$ for $0\leq x\leq2\pi$ and (d) $c_n(3x)$ for $0\leq x\leq\pi$ when $n=20$. 
}}
	\end{center}
\end{figure}
\medskip

\noindent{\bf 2.1\ \ Properties of $s_n(x)$}
\medskip

Some straightforward properties of $s_n(x)$ are:
\medskip

(i)\ \ Since $s_n'(x)=\sum_{k=1}^n \cos kx$, we have the derivatives
\bee\label{e221a}
s_n'(0)=n,\qquad s_n'(\pi)=\sum_{k=1}^n(-)^k=\left\{\begin{array}{ll} 0 & (n\ \mbox{even})\\ -1 & (n\ \mbox{odd}).\end{array}\right.
\ee
It is seen that the curves $s_n(x)$ approach the endpoint $x=\pi$ in different manners according to the parity of $n$.
\medskip

(ii)\ \ With $p$ a positive integer, the functions $s_n(px)$ approach the sawtooth function as $n\to\infty$, with positive and negative parts in $(0, \pi)$ when $p\geq 2$. The extreme values are approximately $\pm 1.852$ as $n\to\infty$ due to the well-known Gibbs phenomenon. In Section 4, we shall take the extreme values to be
\bee\label{e221b}
\xi:=\sup\{\pm s_n(x): x\in (0,2\pi)\}=2. 
\ee
%\medskip

(iii)\ For odd integer $p\geq 3$, the functions $s_n(px)$ are positive in the intervals 
\bee\label{e221c}
(0, \pi/p),\ \ (2\pi/p,\, 3\pi/p),\ \  (4\pi/p,\, 5\pi/p),\,\ldots ,\,((p-1)\pi/p,\,\pi).
\ee

More recondite properties resulting from detailed examination of the curves $s_n(x)$ for different $n$ (which we do not show here) are:
\medskip

(iv)\ For $n\geq 2$, the minimum value of $s_n(x)$ in the interval $[\f{1}{3}\pi, \f{2}{3}\pi]$ is given by
\bee\label{e221d}
s_2(\f{2}{3}\pi)=s_3(\f{2}{3}\pi)=\frac{\sqrt{3}}{4}~.
\ee
%\vspace{0.1cm}

(v)\ For $n\geq 2$, the lowest curve in the intervals $(0,\pi\lambda)$ and $(\f{2}{3}\pi, \pi)$ is 
\bee\label{e221e}
s_2(x)=\sin x(1+\cos x).
\ee
The quantity $\lambda$, which results from the first intersection (measured from $x=0$) of the curves corresponding to $n=2$ and $n=9$, has the value $\lambda\doteq 0.207685$.

\vspace{0.3cm}

\noindent{\bf 2.2\ \ Properties of $c_n(x)$}
\medskip

From the periodicity of $c_n(x)$ we  have
\bee\label{e221}
c_n(2\pi k)=\sigma_n\qquad (k=0, 1, 2, \ldots),
\ee
where
\bee\label{e222}
\sigma_n:=\sum_{k=1}^n\frac{1}{k}=\psi(n+1)+\gamma,
\ee
with $\psi(x)=\g'(x)/\g(x)$ and $\gamma=0.5772156\ldots$ being the Euler-Mascheroni constant. In addition, we have the symmetry property
\bee\label{e223}
c_n(x)=c_n(2\pi-x)\qquad(0\leq x\leq 2\pi).
\ee
The cosine series is then easily seen to satisfy the inequality
\bee\label{e224}
-1<c_n(x)\leq\sigma_n\qquad (x>0).
\ee
Typical plots of $c_n(x)$ and $c_n(3x)$ are shown in Fig.~1.

From the large-$x$ behaviour  $\psi(x)\sim \log\,x-\frac{1}{2x}+O(x^{-2})$
\cite[p.~140]{DLMF}, we see that
\bee\label{e227}
\sigma_n=\log\,n+\gamma+\frac{1}{2n}+O(n^{-2})\qquad (n\to\infty).
\ee
Hence it follows that 
\[c_n(2\pi k)= \log\,n+\gamma+O(n^{-1}) \qquad (k=0, 1, 2, \ldots;\ n\to\infty).\]

Some specific evaluations of $c_n(x)$ will be useful in what follows.
When $x=\pi$ we have the value
\[c_n(\pi)=\sum_{k=1}^n\frac{(-)^k}{k}=\psi(\fs n+1-\fs r)-\psi(n+1),\]
where $r=0$ ($n$ even) and $r=1$ ($n$ odd). 
The value of $c_n(\f{2}{3}\pi)$ is given by
\[c_n(\f{2}{3}\pi)=\sum_{k=1}^n\frac{\cos (\f{2}{3}\pi k)}{k},\]
which, when $n=3m+r$, $r=0, 1, 2$, yields
\begin{eqnarray}
c_n(\f{2}{3}\pi)&=&\frac{1}{3}\bl(1+\frac{1}{2}+\frac{1}{3}+\cdots+\frac{1}{m}\br)-\frac{1}{2}\bl(1+\frac{1}{2}+\frac{1}{4}+\frac{1}{5}+\cdots+\frac{1}{3m+r}\br)\nonumber\\
&=&\frac{1}{2}\bl(1+\frac{1}{2}+\frac{1}{3}+\cdots+\frac{1}{m}\br)-\frac{1}{2}\bl(1+\frac{1}{2}+\frac{1}{3}+\cdots+\frac{1}{3m+r}\br)\nonumber\\
&=&\frac{1}{2}\bl\{\psi(\f{1}{3}n+1-\f{1}{3}r)-\psi(n+1)\br\}.\label{e225}
\end{eqnarray}

Consequently, we find the large-$n$ behaviour
\bee\label{e228}
c_n(\f{2}{3}\pi)= -\frac{1}{2}\log\,3+O(n^{-1}),\qquad c_n(\pi)= -\log\,2+O(n^{-1}) \qquad (n\to\infty).
\ee
\vspace{0.6cm}

\begin{center}
{\bf 3. \ The trigonometric sums $S_n^{(j)}(x)$ ($j=1, 2, 3$)}
\end{center}
\setcounter{section}{3}
\setcounter{equation}{0}
\renewcommand{\theequation}{\arabic{section}.\arabic{equation}}
In Fig.~2 we show typical plots of the trigonometric sums $S_n^{(j)}(x)$ for $1\leq j\leq 4$ and $0\leq x\leq\pi$.
We note that the sums $S_n^{(1)}(x)$ and $S_n^{(2)}(x)$ are symmetrical about $x=\fs\pi$, whereas the other sums do not possess this symmetry. In Fig.~2(a, b) the appearance of a `spike' at $x=\fs\pi$ is just visible when $n=200$; this `spike' becomes more pronounced as $n$ increases. A `spike' is also visible near $x=\f{2}{3}\pi$ in the graph of $S_n^{(3)}(x)$ in Fig.~2(c). The graph of $S_n^{(4)}(x)$, however, exhibits a finite jump near $x=\f{2}{3}\pi$; this is considered in detail in Section 4. All the curves consist of tiny undulations that are too small to be visible on the scale of the figures.

\begin{figure}[th]
	\begin{center}	{\tiny($a$)}\includegraphics[width=0.375\textwidth]{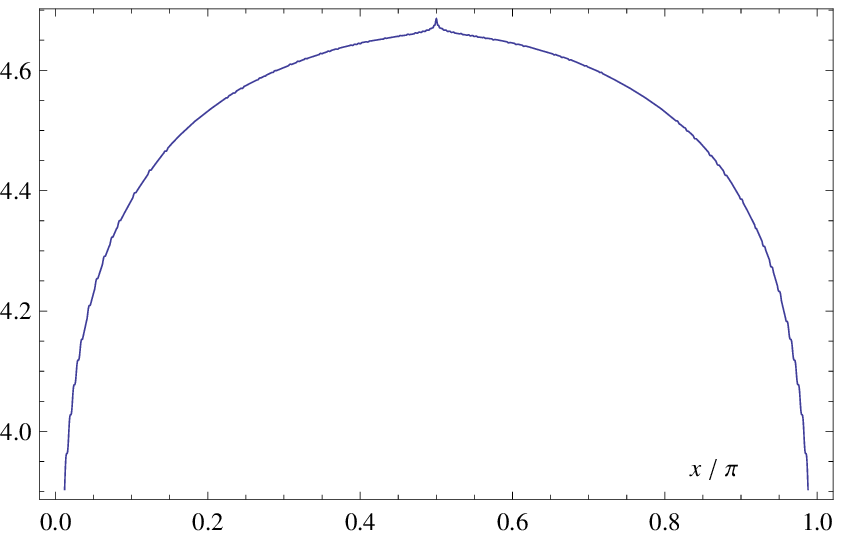}
	\qquad
	{\tiny($b$)}\includegraphics[width=0.375\textwidth]{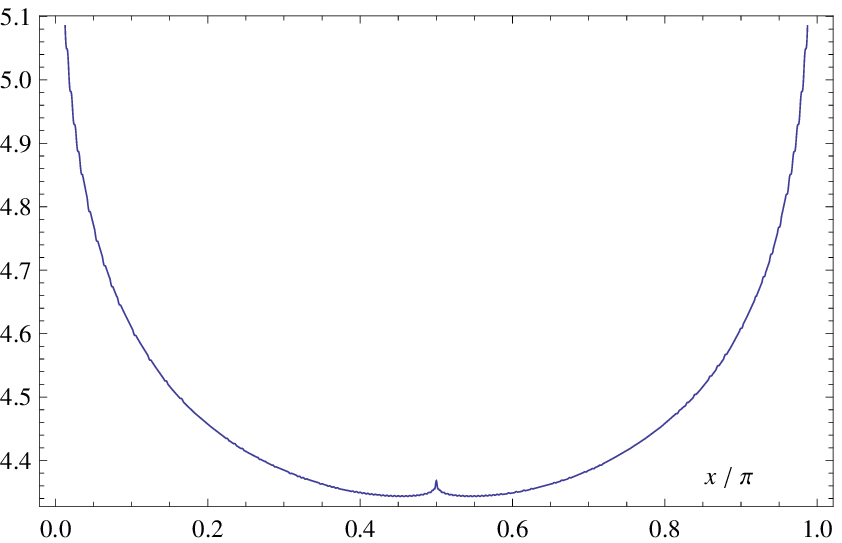}
	\vspace{0.4cm}
	
	{\tiny($c$)}\includegraphics[width=0.375\textwidth]{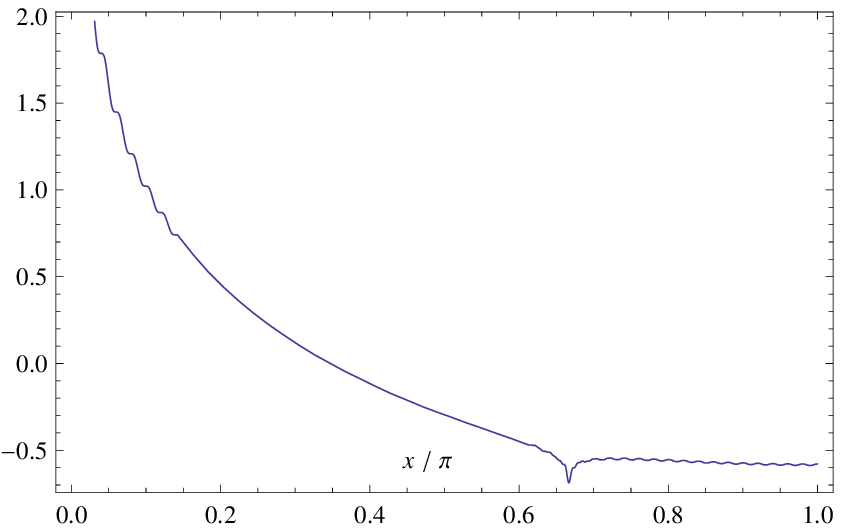}\qquad
	{\tiny($d$)}\includegraphics[width=0.375\textwidth]{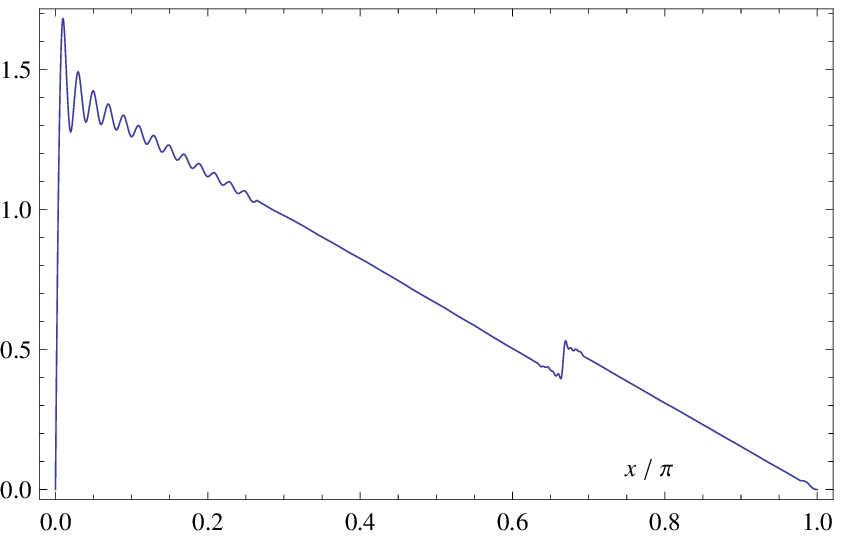}
\caption{\small{Plots of the F\'ejer-type sums for $0\leq x\leq\pi$: (a) $S_n^{(1)}(x)$, (b) $S_n^{(2)}(x)$, (c) $S_n^{(3)}(x)$ when $n=200$ and
(d) $S_n^{(4)}(x)$ when $n=100$. 
}}
	\end{center}
\end{figure}
\vspace{0.3cm}

\noindent {\bf 3.1\ The sum $S_n^{(1)}(x)$}
\vspace{0.3cm}

\noindent The sum $S_n^{(1)}(x)$ satisfies $S_n^{(1)}(0)=\sigma_n \cos (1)$ and, since $\cos kx\in[-1,1]$, it readily follows that
\bee\label{e21}
S_n^{(1)}(x)=\sum_{k=1}^n\frac{\cos (\cos kx)}{k}\geq \cos (1) \sum_{k=1}^n\frac{1}{k}=\sigma_n \cos (1).
\ee
It is found from numerical calculations when $n$ is even and $n\geq18$ that the maximum value of $S_n^{(1)}(x)$ occurs at $x=\fs\pi$; the maximum value also occurs at $x=\fs\pi$ for all odd values of $n$. When $n$ is even and $n\leq16$ the maximum value of 
$S_n^{(1)}(x)$ occurs at the neighbouring local maxima situated at $$x=\frac{1}{2}\pi\pm \frac{\pi}{2n+2}.$$
In Fig.~3(a) we show details of the `spike' occurring in $S_n^{(1)}(x)$ near $x=\fs\pi$ that becomes progreesively more pronounced as $n$ increases.
\begin{figure}[th]
	\begin{center}	{\tiny($a$)}\includegraphics[width=0.375\textwidth]{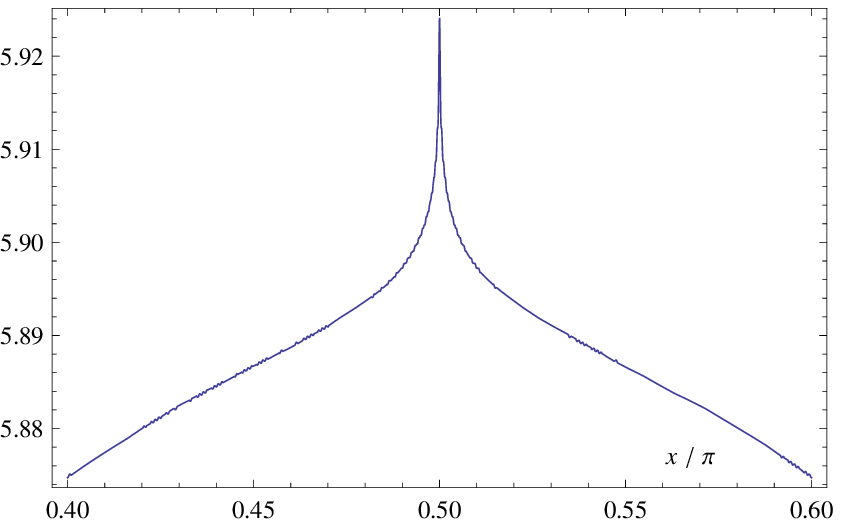}
	\qquad
	{\tiny($b$)}\includegraphics[width=0.375\textwidth]{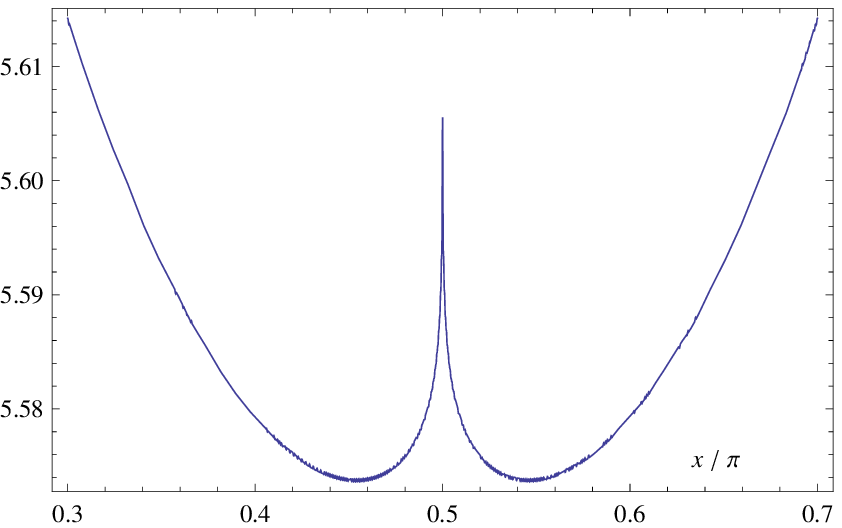}
\caption{\small{Plots of $S_n^{(j)}(x)$, $j=1, 2$ for $x$ near $\fs\pi$ when $n=1000$.
}}
	\end{center}
\end{figure} 

To determine the height of the `spike' we consider the value of $S_n^{(1)}(x)$ at $x=\fs\pi$.
When $n=2m$ is an even integer, we have
\begin{eqnarray*}
S_n^{(1)}(\fs\pi)&=& 1+\frac{1}{3}+\frac{1}{5}+\cdots +\frac{1}{2m-1}+\frac{\cos (1)}{2}\bl(1+\frac{1}{2}+\frac{1}{3}+\cdots +\frac{1}{m}\br)\\
&=& \frac{1}{2}\psi(\fs n+\fs)+\frac{1}{2}\gamma+\log\,2+\frac{\cos (1)}{2}\bl\{\psi(\fs n+1)+\gamma\br\}\qquad (n\ \mbox{even})
\end{eqnarray*}
and in a similar manner when $n$ is an odd integer
\[S_n^{(1)}(\fs\pi)=\frac{1}{2}\psi(\fs n+1)+\frac{1}{2}\gamma+\log\,2+\frac{\cos (1)}{2}\bl\{\psi(\fs n+\fs)+\gamma\br\}\qquad (n\ \mbox{odd}).\]
Thus we have the large-$n$ behaviour
\bee
S_n^{(1)}(\fs\pi)\sim \frac{(1+\cos (1))}{2}\,\log\,(\fs n)\qquad (n\to\infty).
\ee
Thus the height of the `spike' increases without limit as $n\to\infty$ and consequently there is no finite upper bound for the sum $S_n^{(1)}(x)$.
\vspace{0.3cm}

\noindent {\bf 3.2\ The sum $S_n^{(2)}(x)$}
\vspace{0.3cm}

\noindent The sum $S_n^{(2)}(x)$  satisfies $S_n^{(2)}(0)=\sigma_n$ and, since $\sin kx\in[-1,1]$, it readily follows that
\bee\label{e2100}
S_n^{(2)}(x)=\sum_{k=1}^n\frac{\cos (\sin kx)}{k}> \cos (1) \sum_{k=1}^n\frac{1}{k}=\sigma_n \cos (1).
\ee

Plots of $S_n^{(2)}(x)$ are shown in Figs.~2(b) and 3(b), where when $n=200$ the central spike at $x=\fs\pi$ is just visible, whereas when $n=1000$ it is more prominent.  When $n=2m$ is an even integer, the value at $x=\fs\pi$ 
is given by
\begin{eqnarray*}
S_n^{(2)}(\fs\pi)&=&\cos (1)\bl\{1+\frac{1}{3}+\frac{1}{5}+\cdots + \frac{1}{2m-1}\br\}+\frac{1}{2}\bl\{1+\frac{1}{2}+\frac{1}{3}+\cdots +\frac{1}{m}\br\}\\
&=&\cos (1)\bl\{\frac{1}{2}\psi(\fs n+\fs)+\frac{1}{2}\gamma+\log\,2\br\}+\frac{1}{2}\bl\{\psi(\fs n+1)+\gamma\br\}\qquad (n\ \mbox{even})
\end{eqnarray*}
and similarly
\[S_n^{(2)}(\fs\pi)
=\cos (1)\bl\{\frac{1}{2}\psi(\fs n+1)+\frac{1}{2}\gamma+\log\,2\br\}+\frac{1}{2}\bl\{\psi(\fs n+\fs)+\gamma\br\}\qquad (n\ \mbox{odd}).\]

Thus, we have the large-$n$ behaviour of the central spike given by
\bee
S_n^{(2)}(\fs\pi)\sim \frac{(1+\cos (1))}{2}\,\log\,(\fs n)\qquad (n\to\infty),
\ee
which is asymptotically less than the value $S_n^{(2)}(0)=\sigma_n\sim \log\,n$.
This leads us to the conjecture based on  the above estimate of the central spike and numerical evidence that
\bee
\sigma_n \cos (1)<S_n^{(2)}(x)\leq\sigma_n \qquad (0\leq x\leq\pi).
\ee
\vspace{0.3cm}

\noindent {\bf 3.3\ The sum $S_n^{(3)}(x)$}
\vspace{0.3cm}

\noindent The sum $S_n^{(3)}(x)$ satisfies $S_n^{(3)}(0)=\sigma_n \sin (1)$ and
\[S_n^{(3)}(\pi)=\sum_{k=0}^n \frac{\sin (\cos \pi k)}{k}=\sin (1) \sum_{k=1}^n \frac{(-)^k}{k}
=\sin (1)\bl\{\psi(\fs n+1-\fs r)-\psi(n+1)\br\},\]
where $r=0$ ($n$ even) and $r=1$ ($n$ odd). Thus for large $n$ we find
\[S_n^{(3)}(\pi)\sim -\sin (1)\,\log\,2\doteq -0.583263.\]
A typical plot of $S_n^{(3)}(x)$ is shown in Fig.~2(c), where it is seen that there is a negative spike situated near $x=\f{2}{3}\pi$. 

Routine algebra shows that when $n=3m+r$, $m=0, 1, 2, \ldots$ and for $r=0, 1, 2$
\begin{eqnarray*}S_n^{(3)}(\f{2}{3}\pi)&=&\frac{1}{3} \sin (1)\bl(1+\frac{1}{2}+\cdots +\frac{1}{m}\br)-\sin (\fs)\bl(1+\frac{1}{2}+\frac{1}{4}+\frac{1}{5}+\cdots +\frac{1}{3m-1}+\sum_{k=1}^r\frac{1}{3m+k}\br)\\
&=&\frac{1}{3}(\sin(\fs)+\sin (1))\bl(1+\frac{1}{2}+\cdots +\frac{1}{m}\br)-\sin(\fs)\sum_{k=1}^{3m+r}\frac{1}{k}\\
&=&\frac{1}{3}(\sin(\fs)+\sin (1))\bl(\psi(\frac{n}{3}+1-\frac{r}{3})+\gamma\br)-\sin (\fs) \bl(\psi(n+1)+\gamma\br).
\end{eqnarray*}
Thus, we obtain the large-$n$ behaviour
\bee
S_n^{(3)}(\f{2}{3}\pi) \sim \frac{1}{3}\bl\{\sin (1)-2 \sin (\fs)\br\} \log\,n \sim -0.0391267 \log\,n \qquad (n\to\infty).
\ee
This result shows that the negative spike grows without bound as $n\to\infty$. For example, when $n=4\times 10^5$ we find the value $S_n^{(3)}(\f{2}{3}\pi)\doteq -1.011006$.

We now give a qualitative explanation for the appearance of the negative spike in $S_n^{(3)}(x)$ at $x=\f{2}{3}\pi$, which is just visible in Fig.~2(c) when $n=200$. From the first relation in (\ref{e15}) we have
\bee\label{e3B}
S_n^{(3)}(x)=2\sum_{r=0}^\infty (-)^r J_{2r+1}(1)\,c_n((2r+1)x),
\ee
where $c_n(x)$ is defined in (\ref{e12}). The Bessel function satisfies $J_{\nu}(1)>0$ for $\nu\geq 0$ and decreases rapidly with increasing $\nu$ given by the asymptotic formula \cite[(10.19.1)]{DLMF}
\[J_\nu(1)\sim\frac{1}{\sqrt{\pi e}}\bl(\frac{2\nu}{e}\br)^{\!\!-\nu-1/2}\qquad (\nu\to+\infty).\]
The first few values of $J_{2r+1}(1)$ are:
\[J_1(1)\doteq 0.4400510,\quad J_3(1)\doteq 0.0195634,\quad J_5(1)\doteq 2.497577\times 10^{-4},\]
\[ J_7(1)\doteq 1.502326\times 10^{-6},\quad J_9(1)\doteq 5.249250\times 10^{-9}, \ldots\ .\]

Whenever the argument of the $c_n(x)$ function in (\ref{e3B}) equals a multiple of $2\pi$ (where its value is $\sigma_n$) we can expect a `spike' in the graph of $S_n^{(3)}(x)$ as $n\to\infty$. Thus
when $x=\f{2}{3}\pi$, we have $c_n(3x)=c_n(9x)=c_n(15x)= \ldots =\sigma_n$ and, using the periodic and symmetry properties of $c_n(x)$ in (\ref{e223}), we obtain
\[S_n^{(3)}(\f{2}{3}\pi)=-2\sigma_n\bl\{J_3(1)\!-\!J_9(1)\!+\!J_{15}(1)\!-\cdots\br\}+2c_n(\f{2}{3}\pi)\bl\{J_1(1)\!+\!J_5(1)\!-\!J_7(1)\!-\!J_{11}(1)\!+\cdots\br\}.\]
From (\ref{e228}), $c_n(\f{2}{3}\pi)=O(1)$ as $n\to\infty$, so that
\[S_n^{(3)}(\f{2}{3}\pi)=-2\sigma_n\bl\{J_3(1)\!-\!J_9(1)\!+\!J_{15}(1)\!-\cdots\br\}+O(1)\sim
-0.0386302\,\sigma_n\qquad (n\to\infty),\]
which results in an increasingly negative spike as $n\to\infty$.

A similar argument shows that when $x=\f{2}{5}\pi,\ \f{4}{5}\pi$, where $c_n(5x)=c_n(15x)=c_n(25x)=\ldots =\sigma_n$, we have
\begin{eqnarray*}S_n^{(3)}(\f{2}{5}\pi)&=&2\sigma_n\bl\{J_5(1)-J_{15}(1)+J_{25}(1)-\cdots\br\}+O(1)\\
&\sim& +0.000499515\sigma_n\qquad(n\to\infty).
\end{eqnarray*}
This results in the appearance of positive spikes in the neighbourhood of $x=\f{2}{5}\pi,\ \f{4}{5}\pi$ as $n\to\infty$, which have been observed in numerical calculations (which we do not present here).
Other spikes will appear, for example at $x=\f{2}{7}\pi$, $\f{4}{7}\pi$ and $\f{6}{7}\pi$, but require increasingly large values of $n$ for them to be visible.

\vspace{0.6cm}

\begin{center}
{\bf 4. \ The trigonometric sum $S_n^{(4)}(x)$}
\end{center}
\setcounter{section}{4}
\setcounter{equation}{0}
\renewcommand{\theequation}{\arabic{section}.\arabic{equation}}
We now examine the central case of $S_n^{(4)}(x)$ defined by
\[S_n^{(4)}(x)=\sum_{k=1}^n\frac{\sin (\sin kx)}{k}\qquad (0<x<\pi).\]
The aim is to present an argument that strongly supports the conjecture (\ref{e14}) and also explains and quantifies the jump in the graphs of $S_n^{(4)}(x)$ in the vicinity of $x=\f{2}{3}\pi$; see Fig.~2(d). It is evident that $S_1^{(4)}(x)>0$ for $0<x<\pi$ so that we can take $n\geq 2$ throughout.
Application of the second equation in (\ref{e15}) shows that $S_n^{(4)}(x)$ can be expressed as the infinite sum involving unit-argument Bessel functions in the form 
\bee\label{e41}
S_n^{(4)}(x)=2\sum_{r=0}^\infty J_{2r+1}(1)\,s_n((2r+1)x),
%2\bl\{J_1(1) s_n(x)+J_3(1) s_n(3x)+ J_5(1) s_n(5x)+\cdots \br\},
\ee
where the sine series $s_n(x)$ is defined in (\ref{e11}).

Let us define for non-negative integer $m$
\[F_m(n;x)=\sum_{k=0}^m J_{2k+1}(1) s_n((2k\!+\!1)x),\qquad T_m(n;x)=\sum_{k=m+1}^\infty J_{2k+1}(1) s_n((2k\!+\!1)x),\]
so that
\bee\label{e42}
S_n^{(4)}(x)=2\{F_m(n;x)+T_m(n;x)\}\qquad (m=0, 1, 2, \ldots).
\ee
The tail $T_m(n;x)$ of the series consists of positive and negative contributions. This can be majorated by taking the most negative value possible for $s_n(x)$ given by $-\xi=-2$; see (\ref{e221b}). Thus the most negative value of
$T_m(n;x)$ is
\[\xi \!\!\sum_{k=m+1}^\infty J_{2k+1}(1)<\sum_{k=m+1}^\infty \frac{2^{-2k}}{(2k+1)!}\]
upon use of the result \cite[(10.14.4)]{DLMF} $J_\nu(1)<2^{-\nu}/\g(1+\nu)$ for $\nu>-\fs$. 

If we define the bound $B_m$ by
\bee\label{e44bd}
B_m:=\sum_{k=m+1}^\infty \frac{2^{-2k}}{(2k+1)!}=2\sinh \fs-\sum_{k=0}^m \frac{2^{-2k}}{(2k+1)!},
\ee
then we obtain
\bee\label{e43}
\hspace{3cm}S_n^{(4)}(x)>2\{F_m(n;x)-B_m\}\qquad (m=0, 1, 2, \ldots).
\ee
The values of $B_m$ decrease monotonically with $m$ and are given asymptotically by
\bee\label{e43a}
B_m\sim \frac{2^{-2m-2}}{(2m+3)!}\{1+O(m^{-2})\}\qquad (m\to\infty).
\ee
Values of $B_m$ for $1\leq m\leq 10$ are shown in Table 1 from which it is apparent that they decrease very rapidly.

We have the following two lemmas:
\newtheorem{lemma}{Lemma}
\begin{lemma}$\!\!\!.$\ The quantity $B_m$ satisfies the upper bound
\bee\label{e43b}
B_{m-1}<\frac{{\cal H}_m}{2m+1},\qquad {\cal H}_m:=\frac{2^{-2m}}{(2m)!}\,\frac{16m^2}{16m^2-1}\quad (m\geq2).
\ee
Furthermore we have the inequality
\bee\label{e43c}
J_1(1) \sin \frac{\pi}{2m+1} >2B_{m-1}\qquad (m\geq 2).
\ee
\end{lemma}

\noindent{\it Proof.}\ \ From the definition (\ref{e44bd}) we have for $m\geq 2$
\begin{eqnarray*}
B_{m-1}&=&\sum_{k=m}^\infty \frac{2^{-2k}}{(2k+1)!}=\frac{2^{-2m}}{(2m+1)!}\bl\{1+\frac{2^{-2}}{(2m\!+\!2)(2m\!+\!3)}+\frac{2^{-4}}{(2m\!+\!2)\ldots (2m\!+\!5)}+\cdots\br\}\\
&<&\frac{2^{-2m}}{(2m+1)!} \sum_{r=0}^\infty (4m^2)^{-r}=\frac{2^{-2m}}{(2m+1)!}\,\frac{16m^2}{16m^2-1},
\end{eqnarray*}
thereby establishing (\ref{e43b}).

The inequality (\ref{e43c}) can be verified by using the lower bound $\sin \theta\geq 2\theta/\pi$ for $\theta\in[0,\fs\pi]$ and the upper bound for $B_{m-1}$ in (\ref{e43b}). Then (\ref{e43c}) is satisfied if
\[J_1(1)\,\frac{2}{\pi}\,\frac{\pi}{2m+1}>\frac{2{\cal H}_m}{2m+1},\]
which yields $J_1(1)>{\cal H}_m$. Since $J_1(1)\doteq 0.44005$ and ${\cal H}_m\leq {\cal H}_2=\f{1}{378}$ ($m\geq 2$) this is seen to be satisfied.\qquad %$\Box$
\begin{lemma}$\!\!\!.$\ When $m=1$ the sum $F_1(n;x)>0$ in $(0,\pi)$.
\end{lemma}

\noindent{\it Proof.}\ \ With $m=1$ we have
\bee\label{e4j}
F_1(n;x)=J_1(1) s_n(x)+J_3(1) s_n(3x).
\ee
Now, since $s_n(x)>0$ in $(0,\pi)$ by (\ref{e11}) and $s_n(3x)>0$ in $(0,\f{1}{3}\pi)$ and $(\f{2}{3}\pi, \pi)$, it follows that $F_1(n;x)>0$ in $(0,\f{1}{3}\pi)\cup (\f{2}{3}\pi,\pi)$. In the interval $[\f{1}{3}\pi, \f{2}{3}\pi]$, the minimum value of $s_n(x)$ is $\sqrt{3}/4$ by  (\ref{e221d}) and $s_n(3x)\leq0$ with the extreme value taken to be $-\xi$ by (\ref{e221b}). Then
\[F_1(n;x)<\frac{\sqrt{3}}{4} J_1(1)-\xi J_3(1)\doteq 0.151421\qquad (x\in [\f{1}{3}\pi, \f{2}{3}\pi])\]
and hence $F_1(n;x)>0$ in $(0,\pi)$. %$\qquad\Box$
\bigskip

We now proceed with our positivity argument. From (\ref{e43}) when $m=1$, we have
\[S_n^{(4)}(x)>2\{F_1(n;x)-B_1\},\qquad B_1=5.23944\times 10^{-4}.\]
The roots of $F_1(n;x)-B_1=0$ will be close to $x=0$ and $x=\pi$, since $B_1\ll 1$. Let ${\cal L}_\lambda$ denote the intervals $(0, \pi\lambda)\cup (\f{2}{3}\pi, \pi)$, where $\lambda=0.207685$; see \S 2.1(v). Then for $x\in {\cal L}_\lambda$, we have $s_n(3x)>0$ so that from (\ref{e4j}) $F_1(n;x)>J_1(1) s_n(x)=F_0(n;x)$. Furthermore, since the lowest curve of $s_n(x)$ in ${\cal L}_\lambda$ corresponds to $n=2$ by (\ref{e221e}), we have
\[F_1(n;x)>F_0(2;x)=J_1(1) \sin x (1+\cos x)\qquad (x\in {\cal L}_\lambda)\]
and consequently
\[S_n^{(4)}(x)>2\{F_0(2;x)-B_1\}\qquad (x\in {\cal L}_\lambda).\]
The roots of $F_0(2;x)-B_1=0$ are
$x_1^-=1.895\times 10^{-4}\pi$ and $x_1^+=0.957\pi$,
which are situated in ${\cal L}_\lambda$. Hence
\bee\label{e44}
S_n^{(4)}(x)>0 \quad\mbox{for}\quad x\in[x_1^-, x_1^+].
\ee

When $m=2$, we have
\[S_n^{(4)}(x)>2\{F_2(n;x)-B_2\},\qquad B_2=3.111\times 10^{-6}.\]
In $(0,\f{1}{5}\pi)$ and $(\f{4}{5}\pi,\pi)$ both $s_n(3x)$ and $s_n(5x)$ are positive. Hence 
$S_n^{(4)}(x)>2\{F_0(n;x)-B_2\}$ for $x\in (0,\f{1}{5}\pi)\cup (\f{4}{5}\pi,\pi)$. These intervals are contained in
${\cal L}_\lambda$ so that $F_0(n;x)$ can be replaced by $F_0(2;x)$ to yield 
\[S_n^{(4)}(x)>2\{F_0(2;x)-B_2\}\qquad (x\in (0,\f{1}{5}\pi)\cup(\f{4}{5}\pi,\pi)).\]
The roots of $F_0(2;x)-B_2=0$ are $x_2^-=1.125\pi\times 10^{-6}$ and $x_2^+=0.992\pi$. Hence $S_n^{(4)}(x)>0$ in the intervals $[x_2^-, \f{1}{5}\pi)\cup(\f{4}{5}\pi,x_2^+]$. Since $x_2^-<x_1^-<\f{1}{5}\pi$ and $\f{4}{5}\pi<x_1^+<x_2^+$, and $S_n^{(4)}(x)>0$ in $[x_1^-, x_1^+]$ by (\ref{e44}), these intervals overlap and so we deduce that
\bee\label{e45}
S_n^{(4)}(x)>0\quad\mbox{for}\quad x\in [x_2^-, x_2^+].
\ee

This process can be continued. From (\ref{e221c}) with $m\geq 2$, we have $s_n((2r+1)x)>0$ for $1\leq r\leq m$ in the intervals
$$L_m:=(0,\f{\pi}{2m+1})\cup (\f{2m\pi}{2m+1},\pi)\subset {\cal L}_\lambda\qquad (m\geq 2).$$
Then $S_n^{(4)}(x)>0$ for $x\in L_m$. 
It follows that $F_m(n;x)>F_0(n;x)>F_0(2;x)$ for $x\in L_m$, so that from (\ref{e43}) we obtain
\[S_n^{(4)}(x)>2\{F_0(2;x)-B_m\}\qquad (x\in L_m;\ m\geq 2).\]

The roots of $F_0(2;x)-B_m=0$ are $x_m^-$ and $x_m^+$. 
From the leading behaviour of $F_0(2;x)=J_1(1) \sin x (1+\cos x)$ given by
\[F_0(2;x)\sim\left\{\begin{array}{ll} 2xJ_1(1) & (x\to0)\\ \fs(\pi-x)^3 J_1(1) & (x\to\pi),\end{array}\right.\]
we find that
\[x_m^-\sim\frac{B_m}{2J_1(1)},\qquad x_m^+\sim \pi- \bl(\frac{2B_m}{J_1(1)}\br)^{\!1/3}\quad (m\gg 1).\]

\begin{table}[t]
\caption{\footnotesize{Values of the roots of $F_0(2,x)-B_m=0$ for $1\leq m\leq 10$.}}
\begin{center}
\begin{tabular}{c|l|ll}
\hline
&&&\\[-0.30cm]

\mcol{1}{c|}{$m$} & \mcol{1}{c|}{$B_m$} & \mcol{1}{c}{$x_m^-/\pi$} & \mcol{1}{c}{$x_m^+/\pi$}\\
[.1cm]\hline
&&&\\[-0.30cm]
1 & $5.2394432083\times10^{-4}$ & $1.8949717051\times10^{-4}$ & 0.9574201720\\
2 & $3.1109874947\times10^{-6}$ & $1.1251639101\times10^{-6}$ & 0.9923025565\\
3 & $1.0789082025\times10^{-8}$ & $3.9021325986\times10^{-9}$ & 0.9988349231\\
4 & $2.4504202963\times10^{-11}$& $8.8625379772\times10^{-12}$& 0.9998468538\\
5 & $3.9253367506\times10^{-14}$& $1.4196930249\times10^{-14}$& 0.9999820808\\
6 & $4.6717514135\times10^{-17}$& $1.6896519501\times10^{-17}$& 0.9999981010\\
8 & $3.1377910235\times10^{-23}$& $1.1348580548\times10^{-23}$& 0.9999999833\\
10& $9.2262804088\times10^{-30}$& $3.3369075760\times10^{-30}$& 0.9999999999\\
[.15cm]\hline
\end{tabular}
\end{center}
\end{table}

We now establish a lemma concerning the roots $x_m^\pm$ and $x_{m-1}^\pm$.
\begin{lemma}$\!\!\!.$\ The roots $x_m^\pm$ and $x_{m-1}^\pm$ satisfy the inequalities
\bee\label{eLem2}
x_m^-<x_{m-1}^-<\frac{\pi}{2m+1},\qquad \frac{2m\pi}{2m+1}<x_{m-1}^+<x_m^+\qquad(m\geq 2).
\ee
\end{lemma}

\noindent{\it Proof.}\ \ Since, from (\ref{e43a}), $B_m$ is a rapidly decreasing function of $m$ it follows that $x_m^-$ and $x_m^+$ are decreasing and increasing sequences, respectively. Values of $x_m^\pm$ and $B_m$ for $1\leq m\leq10$ are presented in Table 1.

The roots $x_m^-$ are determined from $\sin x(1+\cos x)=B_m/J_1(1)$, from which it follows for $x_m^-\in(0,\f{1}{3}\pi)$ that
\[\frac{B_m}{2J_1(1)}<\sin x_m^- <\frac{2B_m}{3J_1(1)}.\]
The requirement that $x_{m-1}^-<\pi/(2m+1)$ is satisfied if $\arcsin (2B_{m-1}/(3J_1(1))<\pi/(2m+1)$, or
\[J_1(1) \sin \frac{\pi}{2m+1}>\frac{2}{3}B_{m-1}\qquad (m\geq 2).\]
This is seen to be verified by (\ref{e43c}) in Lemma 1.

For the roots $x_m^+$, we let ${\hat x}:=\pi-x$. Then $\sin {\hat x}(1-\cos {\hat x})=B_m/J_1(1)$, whence for ${\hat x}_m^+\in (0,\f{1}{3}\pi)$ we have
\[\frac{B_m}{J_1(1)}<\sin {\hat x}_m^+ <\frac{2B_m}{J_1(1)}.\]
The second inequality in (\ref{eLem2}) $x_{m-1}^+>2m\pi/(2m+1)$ is therefore satisfied if $\pi-\arcsin (2B_{m-1}/J_1(1)) >2m\pi/(2m+1)$, or
\[J_1(1) \sin \frac{\pi}{2m+1}>2B_{m-1}\qquad (m\geq 2),\]
which is true by Lemma 1.
\medskip

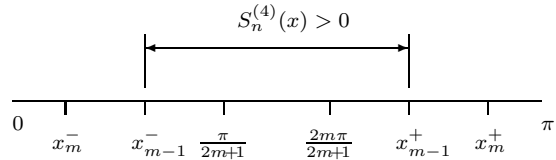
\begin{figure}[h]
\centering
\begin{picture}(200,90)(0,0)
\put(0,40){\line(1,0){200}}
\put(50,60){\line(1,0){100}}
\put(55,60){\vector(-1,0){5}}
\put(145,60){\vector(1,0){5}}
\put(50,45){\line(0,1){20}}
\put(150,45){\line(0,1){20}}
\put(85,68){\footnotesize{$S_n^{(4)}(x)>0$}}
\put(50,40){\line(0,-1){5}}
\put(150,40){\line(0,-1){5}}
\put(180,40){\line(0,-1){5}}
\put(20,40){\line(0,-1){5}}
\put(80,40){\line(0,-1){5}}
\put(120,40){\line(0,-1){5}}
\put(45,22){\footnotesize{$x_{m-1}^-$}}
\put(145,22){\footnotesize{$x_{m-1}^+$}}
\put(15,22){\footnotesize{$x_m^-$}}
\put(175,22){\footnotesize{$x_m^+$}}
\put(70,22){\footnotesize{$\frac{\pi}{2m\!+\!1}$}}
\put(110,22){\footnotesize{$\frac{2m\pi}{2m\!+\!1}$}}
\put(0,29){\footnotesize{$0$}}
\put(200,29){\footnotesize{$\pi$}}
\end{picture}
\caption{\small{Schematic representation of the roots $x_m^\pm$ and $x_{m-1}^\pm$ and the interval of positivity $[x_{m-1}^-, x_{m-1}^+]$.}}
\end{figure}

Since the intervals $[x_m^-,\pi/(2m\!+\!1))$ and $(2m\pi/(2m\!+\!1), x_m^+]$ overlap with the interval of positivity $[x_{m-1}^-,x_{m-1}^+]$  from the previous step (see Fig.~4), it then follows that 
\[S_n^{(4)}(x)>0 \qquad  x\in [x_m^-,\ x_m^+]\]
thereby establishing positivity for the $m$th step.

From Table 1, it is seen that 
when $m=10$, we have $x_{10}^-=3.337\pi\times 10^{-30}$ and $x_{10}^+=0.999\,999\,999\,9\pi$. Thus, with this positivity argument we can approach the endpoints $x=0$ and $x=\pi$ as close as we please.
\vspace{0.6cm}
\newpage

\begin{center}
{\bf 5. \ Concluding remarks}
\end{center}
\setcounter{section}{5}
\setcounter{equation}{0}
\renewcommand{\theequation}{\arabic{section}.\arabic{equation}}
We explain and quantify the positive jump in the graph of $S_n^{(4)}(x)$ appearing in the vicinity of $x=\f{2}{3}\pi$ that is visible in Fig.~2(d). From (\ref{e41}), the functions $s_n(3x)$, $s_n(9x)$, $s_n(15x), \ldots$ have arguments close to $2\pi$, $6\pi$, $10\pi, \ldots\ $, where they undergo a rapid change from $-\Delta$ to $+\Delta$ (the Gibbs phenomenon), where $\Delta\doteq1.852$ as $n\to\infty$; see \S 2.1(ii). Thus we expect a positive jump in the neighbourhood of $x=\f{2}{3}\pi$ equal to approximately
\bee\label{e52}
2\{J_3(1)+J_9(1)+J_{15}(1)+\cdots\} 2\Delta \doteq 0.1449\qquad (n\to\infty).
\ee
Numerical calculations with $n$ large confirm this value for the jump.

Finally, 
since $d/dx\, S_n^{(4)}(x)=\sum_{k=1}^n \cos (\sin kx) \cos kx$, we have the derivatives
\bee\label{e51}
d/dx\,S_n^{(4)}(0)=n,\qquad d/dx\,S_n^{(4)}(\pi)=\left\{\begin{array}{ll} 0 & (n\ \mbox{even})\\ -1 & (n\ \mbox{odd}).\end{array}\right.
\ee
Thus the gradient of $S_n^{(4)}(x)$ at $x=0$ becomes progressively steeper as $n\to\infty$ and the gradient at $x=\pi$ is either 0 or $-1$ depending on the parity of $n$.
The expansions of $\sin (\sin kx)$ near $x=0$ and $x=\pi$ are
\[\sin (\sin kx)=\left\{\begin{array}{ll} kx-\f{1}{3}(kx)^3+\f{1}{10}(kx)^5-\cdots & (x\to 0) \\
\\
(-1)^{k-1}(k{\hat x}-\f{1}{3}(k{\hat x})^3+\f{1}{10}(k{\hat x})^5-\cdots ) & (x\to\pi),\end{array}\right.\]
where ${\hat x}:=\pi-x$. This leads to the expansions
\[S_n^{(4)}(x)=nx-\f{1}{6}n(n+1)(2n+1)\bl(\f{1}{3}x^3-\f{1}{50}(3n^2+3n-1)x^5+\cdots\br) \quad (x\to0)\]
and
\[S_n^{(4)}(x)=\left\{\!\!\begin{array}{l}0\\ {\hat x}\end{array}\!\!\right\}+(-)^n n(n+1)\bl(\f{1}{6}{\hat x}^3-\f{1}{20}(n^2+n-1) {\hat x}^5+\cdots\br)
\qquad n\left\{\begin{array}{l}\!\!\mbox{even}\\\!\!\mbox{odd}\end{array}\right. \quad (x\to\pi),\]
These results are readily seen to yield the derivatives in (\ref{e51}).


\vspace{0.6cm}

\begin{thebibliography}{99}
\footnotesize{
\bibitem{Alz}
H. Alzer, Private communication, 2019.

\bibitem{Fej}
L. F\'ejer, Einige S\"{a}tze, die sich auf das Vorzeichen einer ganzen rationalen Funktion beziehen. 
%nebst Anwendungen dieser S\"{a}tze auf die Abschnitte und Abschnittsmittelwerte von ebenen und r\"{a}umlichen harmonischen Entwicklungen und von beschr\"{a}nkten Potenzriehen. 
Monatsh. Math. {\bf 35} (1928) 305--344.

\bibitem{Gron}
T.H. Gronwall, \"{U}ber die Gibbssche Erscheinung und die trigonometrische Summe $\sin x+\fs \sin 2x+\cdots +\f{1}{n} \sin nx$. Math. Anal. {\bf 72} (1912) 228--243.

\bibitem{Jack}
D. Jackson, \"{U}ber die trigonometrische Summe, Rend. Circ. Math. Palermo {\bf 32} (1911) 257--262.

\bibitem{Lan}
E. Landau, \"{U}ber eine trigonometrische Ungleichung. Math. Z. {\bf 37} (1933) 36.

\bibitem{Tur}
P. Tur\'an, On a trigonometric sum. Ann. Soc. Polon. Math. {\bf 25} (1952) 155-161.

\bibitem{DLMF}
F.W.J. Olver, D.W. Lozier, R.F. Boisvert, C.W. Clark (eds.), {\it   
NIST Handbook of Mathematical Functions}, Cambridge University Press, Cambridge, 2010.

%\bibitem{PK}
%R.B. Paris, D. Kaminski, {\it Asymptotics and Mellin-Barnes Integrals}, Cambridge University Press, Cambridge, 2011.

\bibitem{WHY}
W.H. Young, On certain series of Fourier, Proc. London Math. Soc. {\bf 11} (1913) 357--366.

%\bibitem{WBF}
%G.N. Watson, {\it A Treatise on the Theory of Bessel Functions}, Cambridge University Press, Cambridge 1952.

}
\end{thebibliography}
\end{document}